\documentclass[11pt]{article}
\usepackage{epic,latexsym,amssymb}
\usepackage{color}
\usepackage{tikz}

\usepackage{geometry,graphicx,verbatim,amsmath}
\usepackage{tikz,ifthen}
\usetikzlibrary{calc}

\newtheorem{thm}{Theorem}[section]
\newtheorem{lem}[thm]{Lemma}
\newtheorem{cor}[thm]{Corollary}

\newtheorem{prop}[thm]{Proposition}

\newcommand{\rt}{{\rm rt}}
\newcommand{\mdim}{{\rm mdim}}
\newcommand{\edim}{{\rm edim}}

\newcommand{\proof}{\noindent{\bf Proof.\ }}
\newcommand{\qed}{\hfill $\square$ \bigskip}

\let\oldenumerate\enumerate
\renewcommand{\enumerate}{
  \oldenumerate
  \setlength{\itemsep}{1pt}
  \setlength{\parskip}{0pt}
  \setlength{\parsep}{0pt}
}

\begin{document}

\title{On mixed metric dimension in subdivision, middle, and total graphs}

\author{$^{a}$Ali Ghalavand
\and $^{b,c,d}$Sandi Klav\v zar
\and $^{a}$Mostafa Tavakoli
\and $^e$Ismael G. Yero
\\ \\
$^a$Department of Applied Mathematics, Faculty of Mathematical Sciences,\\
Ferdowsi University of Mashhad, Iran \\
\small \tt alighalavand@grad.kashanu.ac.ir, \small \tt m$\_$tavakoli@um.ac.ir\\
$^b$ Faculty of Mathematics and Physics, University of Ljubljana, Slovenia\\
$^c$ Faculty of Natural Sciences and Mathematics, University of Maribor, Slovenia\\
$^d$ Institute of Mathematics, Physics and Mechanics, Ljubljana, Slovenia\\
\small \tt sandi.klavzar@fmf.uni-lj.si\\
$^e$ Departamento de Matem\'{a}ticas, Universidad de C\'{a}diz, Algeciras, Spain\\
\small \tt ismael.gonzalez@uca.es\\
}

\date{}
\maketitle

\begin{abstract}
Let $G$ be a graph and let $S(G)$, $M(G)$, and $T(G)$ be the subdivision, the middle, and the total graph of $G$, respectively. Let $\dim(G)$, $\edim(G)$, and $\mdim(G)$ be the metric dimension, the edge metric dimension, and the mixed metric dimension of $G$, respectively. In this paper, for the subdivision graph it is proved that $\frac{1}{2}\max\{\dim(G),\edim(G)\}\leq\mdim(S(G))\leq\mdim(G)$. A family of graphs $G_n$ is constructed for which $\mdim(G_n)-\mdim(S(G_n))\ge 2$ holds and this shows that the inequality  $\mdim(S(G))\leq\mdim(G)$ can be strict, while for a cactus graph $G$, $\mdim(S(G))=\mdim(G)$.  For the middle graph it is proved that $\dim(M(G))\leq\mdim(G)$ holds, and if $G$ is tree with $n_1(G)$ leaves, then $\dim(M(G))=\mdim(G)=n_1(G)$. Moreover, for the total graph it is proved that $\mdim(T(G))=2n_1(G)$ and $\dim(G)\leq\dim(T(G))\leq n_1(G)$ hold when $G$ is a tree.
\end{abstract}

{\textbf{AMS Subj.\ Class.:} 05C12; 05C69; 05C76}

\section{Introduction}

The metric dimension of graphs is a classical topic in graph theory that has significantly increased the attention of several researchers in the last decade. This can be seen for instance while making a search query at the MathSciNet database, where about 90 percent of articles were published in the last 10 years. Reasons for this increasing interest might arise from the several applications of metric dimension topics on practical problems from areas like computer science \cite{Melter-1984}, social networks \cite{Trujillo-Rasua-2016}, chemistry \cite{Johnson-1993}, and biology \cite{till-2019}, among other ones. For more information on these facts we suggest the recent survey \cite{till-2022+}. Another reason for this might be regarding the explosion of different newly defined variations of the classical concept that give more insight into the classical topic and between themselves. For more information and background on many of these variants, we recommend the other recent survey \cite{dorota-2022+}.

It is our goal to present some contributions on one of these variants called mixed metric dimension. Particularly the exposition centers the attention on finding (or bounding) the mixed metric dimension of some graphs obtained under some modifications on their structure which have the common nature of making subdivisions into their edges. The modifications are as follows.

\begin{itemize}
  \item The {\em subdivision graph} $S(G)$ of a graph $G$ is obtained from $G$ by subdividing each edge of $G$.  It is clear that $n(S(G)) = n(G) + m(G)$ and $m(S) = 2m(G)$. If $e \in E(G)$, then by $v_{e}$ we will denote the vertex of $S(G)$ obtained by subdividing the edge $e$, and call it a {\em subdivision vertex}.
  \item The {\em middle graph} $M(G)$ of a graph $G$ is obtained from $G$ by subdividing each edge, and then joining pairs of these new vertices if and only if their corresponding edges have a common endvertex in $G$. Note that $M(G)$ contains the line graph $L(G)$ of $G$ as a proper induced subgraph and the subdivision graph $S(G)$ as a proper spanning subgraph. Middle graphs were introduced by Hamada and Yoshimura~\cite{hamada-1976}, and further studied in several papers, cf.~\cite{kwak-2006, somodi-2017}.
  \item The {\em total graph} $T(G)$ of a graph $G$ is obtained from the middle graph $M(G)$ by adding the edges $xy$, where $xy\in E(G)$. Total graphs were introduced more than 50 years ago, see~\cite{behzad-1966}, and further on studied in about two hundred articles. Note that $T(G)$ contains $G$ as an induced subgraph, and its edges will be called {\em original edges}. $T(G)$ also contains $S(G)$ as a spanning subgraph, its edges will be addressed as {\em $S(G)$-edges}. In addition, $T(G)$ contains $L(G)$ as an induced subgraphs, and its edges shall be called {\em $L(G)$-edges}.
\end{itemize}

The subdivision, middle and total graphs are combinatorial constructions whose own properties are many times inherited from the original graphs from which they have arisen. Consequently, these constructions, as well as several other related ones (line graphs, edge contracted graphs, etc), are frequently considered in the investigation on graph theory as they are giving some extra information on the original graphs. The case of metric dimension related topics does not escape to this, as we shall present throughout our whole exposition.

\section{Preliminaries}

In this paper, we consider finite, simple and connected graphs. The vertex and edge sets of a graph $G$ are denoted by $V(G)$ and $E(G)$, respectively. Two vertices $x,y\in V(G)$ are \emph{resolved} (resp.\ \emph{determined} or \emph{identified}) by a vertex $v\in V(G)$ if $d_G(x,v)\ne d_G(y,v)$ where $d_G$ stands for the classical distance operator in graphs. A set $S\subset V(G)$ is called a \emph{resolving set} for $G$ if every two vertices $u,v\in V(G)$ are resolved by a vertex of $S$. A resolving set of the smallest possible cardinality in $G$ is called a \emph{metric basis}, and the cardinality of a metric basis is the metric dimension of $G$, denoted $\dim(G)$. These concepts were first formally and independently introduced for graphs in \cite{HM-1976, S-1975} in connection with some problems of uniquely recognizing the vertices of a graph with a purpose of locating intruders in a network. We may remark that in \cite{S-1975}, resolving sets were named \emph{locating sets} and the metric dimension was called \emph{locating number}.

On the other hand, and aimed to also uniquely identify the edges (resp.\ edges and vertices) of a graph, the notion of edge (resp.\ mixed)  resolving sets were first presented in \cite{KTY-2018} (resp.\ \cite{KKTY-2017}) as follows. Two elements $x,y\in V(G)\cup E(G)$ are \emph{resolved} by a vertex $v\in V(G)$ if $d_G(x,v)\ne d_G(y,v)$. Here, if $x=ww'$ is an edge, then $d_G(x,v)=\min\{d_G(w,v), d_G(w',v)\}$. A set of vertices $S\subset V(G)$ is said to be an \emph{edge resolving set} for $G$ if every two edges $e,f\in E(G)$ are resolved by a vertex of $S$. Moreover, the set $S$ is called a \emph{mixed resolving set} for $G$ if any two elements (vertices or edges) $x,y\in V(G)\cup E(G)$ of the graph are resolved by a vertex of $S$. An edge (resp.\ mixed) resolving set of the smallest possible cardinality is called an edge (resp.\ mixed) metric basis, and its cardinality the edge (resp.\ mixed) metric dimension, denoted by $\edim(G)$ (resp.\ $\mdim(G)$). Recent studies on the edge and mixed metric dimensions of graphs are for instance \cite{FKK-2019,sedlar-2022} and \cite{MM2, MKSM-2021,SS-2021}, respectively.

The following concept from the seminal paper~\cite{KKTY-2017} is very useful. A vertex $ u \in N_G(v)$ is a {\em maximal neighbor of} $v$ if all neighbors of $v$ (and $v$ itself) are in $N_G[u]$. That is, $u$ is a maximal neighbor of $v$ if $u$ is adjacent to all neighbors of $v$. The next Lemma~\ref{lem:max-neigbor} implicitly follows from~\cite[Theorem 3.8]{KKTY-2017}.

\begin{lem}\label{lem:max-neigbor}
If $W$ is a mixed resolving set for a graph $G$, and $v\in V(G)$ has a maximal neighbour, then $v\in W$.
\end{lem}

\proof
Since $v$ has a maximal neighbour, there exists $u\in N_G(v)$ such that $N_G[v]\subseteq N_G[u]$. Consider the edge $uv$ and the vertex $u$. Since for $x\in V(G)\backslash\{v\}$ we have $d_G(x,u) = d_G(x,uv)$, the only vertex that can distinguish $uv$ and $u$ is $v$. That is, we conclude that $v\in W$.
\qed

Recall that a graph $G$ with edge disjoint cycles is called a {\em cactus}. For three vertices $u$, $v$ and $w$ from a cycle $C$ in a cactus graph $G$ we say that they form a {\em geodesic triple} of vertices if
$d_G(u, v) + d_G(v, w) + d_G(w, u) = |V (C)|$. Let $G$ be a cactus graph with $c$ cycles $C_1,\ldots,C_c$. The {\em root vertices} on the cycle $C_i$ are the vertices of degree at least $3$.  The number of root vertices on the cycle $C_i$ is denoted by $\rt(C_i)$. By $n_1$ we denote the number of degree $1$ vertices of $G$.

\begin{thm}\label{th0}{\rm\cite[Theorem 5]{SS-2021}}
Let $G$ be a cactus graph with $t$ cycles $C_{s_1},\ldots,C_{s_t}$. Then
$$\mdim(G)=n_1(G)+\sum_{i=1}^t\max\{3-\rt(C_{s_i}),0\}+\epsilon,$$
where $\epsilon$ is the number of cycles $C_{s_i}$ in $G$ for which $\rt(C_{s_i})\geq3$ and there is not a geodesic triple of root vertices on the cycle $C_{s_i}$.
\end{thm}

\section{Subdivision graphs}
\label{sec:main}

Let $G$ be a graph and let $X$ be a metric basis of $S(G)$. Hence, denote
$$\phi(X) = \{ u\in V(G):\ u\in X\ {\rm or}\ u\ {\rm is\ adjacent\ to\ a\ subdivision\ vertex\ from}\ X\}$$
and
$$\phi(G) = \min \{ |\phi(X)|:\ X\ {\rm is\ a\ metric\ basis\ of}\ S(G)\}\,.$$

\begin{thm}
\label{thm:subdivided}
If $G$ is a graph, then the following hold:
\begin{enumerate}
\item[(i)] $\mdim(S(G)) \le \mdim(G)$,
\item[(ii)] $\phi(G) \ge \max\{ \dim(G), \edim(G)\}$.
\end{enumerate}
\end{thm}

\proof
At the beginning of the proof, we give the distance function for those vertex to vertex and vertex to edge situations in the graph $S(G)$ needed later on. Let $x,y\in V(G)$. Then $xv_1v_2\ldots v_ky$ is a path of length $k+1$ in $G$ if and only if
$xv_{xv_1}v_1v_{v_1v_2}v_2\ldots v_kv_{v_ky}y$ is a path of length $2(k+1)$ in $S(G)$. Thus, if $x,y\in V(G)$, then
\begin{equation}
\label{eq:distance-vertex-vertex}
d_{S(G)}(x,y)=2d_G(x,y)\,.
\end{equation}
Similarly, if $x\in V(G)$ and $e\in E(G)$, then
\begin{equation}
\label{eq:distance-vertex-subidivided-vertex}
d_{S(G)}(x,v_{e}) = 2d_G(x,e)+1\,,
\end{equation}
and if $e, e'\in E(G)$, then
\begin{equation}
\label{eq:distance-subidivided-vertex-subidivided-vertex}
d_{S(G)}(v_{e}, v_{e'}) = 2d_G(e,e')+2\,.
\end{equation}
Finally, if $x\in V(G)$ and $f\in E(S(G))$, where $f$ is one of the edges obtained by subdividing $e\in E(G)$, then
\begin{equation}
\label{eq:distance-vertex-edge}
d_{S(G)}(x, f) \in \{2d_G(x,e),2d_G(x,e)+1\}\,.
\end{equation}

Let us now prove the two claims of the theorem.

(i) Let $W\subseteq V(G)$ be a mixed resolving set for $G$. We are going to show that $W\subseteq V(S(G))$ is a mixed resolving set for $S(G)$.

Let $x,y\in V(S(G))$. If $x,y\in V(G)$, then there exists $w\in W$ such that $d_G(x,w)\ne d_G(y,w)$. By~\eqref{eq:distance-vertex-vertex} we have $d_{S(G)}(x,w) = 2 d_G(x,w) \ne 2d_G(y,w) = d_{S(G)}(y,w)$. Suppose next that $x\in V(G)$ and $y = v_e$ for some $e\in E(G)$. Then by~\eqref{eq:distance-vertex-vertex} and~\eqref{eq:distance-vertex-subidivided-vertex} we get that  $d_{S(G)}(x,w)$ is even and $d_{S(G)}(v_e,w)$ is odd for any $w\in W$, hence $x$ and $v_e$ are distinguished by (every vertex of) $W$. Suppose finally that $x=v_e$ and $y = v_{e'}$ for two edges $e, e'\in E(G)$. Since $W$ is a mixed resolving set, there exists $w\in W$ such that $d_G(e,w)\ne d_G(e',w)$. Using~\eqref{eq:distance-vertex-subidivided-vertex} once more, we get $d_{S(G)}(w,v_e) = 2 d_G(w,e) + 1 \ne 2d_G(w,e') + 1 = d_{S(G)}(w,v_{e'})$. So any pair of vertices of $S(G)$ is distinguished by $W$.

Consider next a vertex $x\in V(S(G))$ and an edge $f\in E(S(G))$. Then $f = zv_e$, where $e=zz'\in E(G)$.  Suppose that $x\in V(G)$ and let $w\in W$ be a vertex with $k = d_G(x,w)\ne d_G(e,w)= \ell$. Then $d_{S(G)}(x,w) = 2k$ and $d_{S(G)}(f,w)\in \{2\ell, 2\ell + 1\}$ by~\eqref{eq:distance-vertex-edge}. Hence $d_{S(G)}(x,w) \ne d_{S(G)}(f,w)$. Suppose next $x\in V(S(G))\setminus V(G)$. Then $x = v_{e'}$, where $e' \in E(G)$. As $W$ is a mixed resolving set, there exists $w\in W$ such that $k = d_G(e,w) \ne d_G(e',w)= \ell$. Then $d_{S(G)}(x,w) = 2\ell+1$ by~\eqref{eq:distance-vertex-subidivided-vertex} and $d_{S(G)}(f,w)\in \{2k, 2k + 1\}$ by~\eqref{eq:distance-vertex-edge}. Hence again $d_{S(G)}(x,w) \ne d_{S(G)}(f,w)$.

It remains to verify that also each pair of edges of $S(G)$ is distinguished by $W$. Let $f = zv_e\in E(S(G))$, where $e=zz'\in E(G)$ and $f' = yv_{e'}\in E(S(G))$, where $e'=yy'\in E(G)$. Let $w\in W$ be a vertex with $k = d_G(w,e)\ne d_G(w,e')= \ell$. Then by~\eqref{eq:distance-vertex-edge}, $d_{S(G)}(w,f)\in \{2k, 2k + 1\}$ and $d_{S(G)}(w,f')\in \{2\ell, 2\ell + 1\}$ and consequently  $d_{S(G)}(f,w) \ne d_{S(G)}(f',w)$. This proves the first inequality.

\medskip
(ii) In order to establish $\phi(G) \ge \dim(G)$, we need to prove that if $X$ is an arbitrary metric basis of $S(G)$, then $\phi(X)$ is a resolving set of $G$. Consider $u, v\in V(G)$. If at least one of them belongs to $X$, then it also belongs to $\phi(X)$; so $u$ and $v$ are distinguished by it.  Suppose this is not the case. Then in $S(G)$ there exists a vertex $x\in X$ such that $d_{S(G)}(u,x) \ne d_{S(G)}(v,x)$. If $x\in V(G)$, then $x\in \phi(X)$ and by~\eqref{eq:distance-vertex-vertex} we are done. Suppose $x$ is a subdivision vertex, say $x = v_e$, where $e = x'x''$. Then $x', x''\in \phi(X)$. If $d_G(u,x') = d_G(v,x')$ and $d_G(u,x'') = d_G(v,x'')$, then by~\eqref{eq:distance-vertex-edge} we would have $d_{S(G)}(u,x) = d_{S(G)}(v,x)$. Hence $u, v$ are distinguished in $G$ by at least one of $x', x''$, a contradiction. Thus $\phi(G) \ge \dim(G)$ holds.

It remains to prove that $\phi(G) \ge \edim(G)$. Consider $e=xx', f=yy'\in E(G)$. Then consider the vertices $v_e$ and $v_f$ of $S(G)$. If at least one of these two vertices lies in $X$, say $v_e$, then $d_G(e,x) = d_G(e,x') = 0$, and $d_G(f,x) \ne 0$ or $d_G(f,x') \ne 0$, hence $e$ and $f$ are distinguished by $\phi(X)$ since $x,x'\in \phi(X)$. It remains to consider the case when $v_e\notin X$ and $v_f\notin X$. Then in $S(G)$ there exists a vertex $w\in X$ such that $d_{S(G)}(w,v_e) \ne d_{S(G)}(w,v_f)$. If $w\in V(G)$, then $w\in \phi(X)$ and by~\eqref{eq:distance-vertex-subidivided-vertex} we are done. Suppose now that $w$ is a subdivision vertex, say $w = v_g$, where $g = w'w''$. Then $w', w''\in \phi(X)$. By~\eqref{eq:distance-subidivided-vertex-subidivided-vertex} we have $d_{S(G)}(v_g,v_e) = 2k+2$ (for some $k$) and $d_{S(G)}(v_g,v_f) = 2\ell+2$ (for some $\ell$), where $k\ne \ell$. We may assume without loss of generality that $\ell < k$. Again assume without loss of generality that $d_G(w'',f) \le d_G(w',f)$. Then $d_G(g,f) = d_G(w'',f) = \ell$. On the other hand, because $d_G(g,e) = k$, we have $d_G(w'',e) \ge k > \ell = d_G(w'',f)$. We conclude that $w''\in \phi(X)$ distinguishes the edges $e$ and $f$.
\qed

\begin{cor}
\label{cor:subdivided}
If $G$ is a graph, then
$$\frac{1}{2}\max \{\dim(G), \edim(G)\} \le \frac{1}{2}\phi(G) \le \mdim(S(G))\le \mdim(G)\,.$$
\end{cor}

\proof 
Let $X$ be an arbitrary metric basis of $S(G)$. Then by the construction,
\begin{equation}\label{eq:cor1}
\phi(G) \le |\phi(X)| \le 2|X| = 2\dim(S(G))\,.
\end{equation}
Then we have
$$\frac{1}{2}\max \{\dim(G), \edim(G)\} \le \frac{1}{2}\phi(G)\le \dim(S(G)) \le \mdim(S(G))\le \mdim(G)\,$$
where the first and the last inequality hold by Theorem~\ref{thm:subdivided} (ii) and (i), respectively. The second inequality was established in \eqref{eq:cor1}, and the third inequality is obvious by the definition of the two dimensions.
\qed

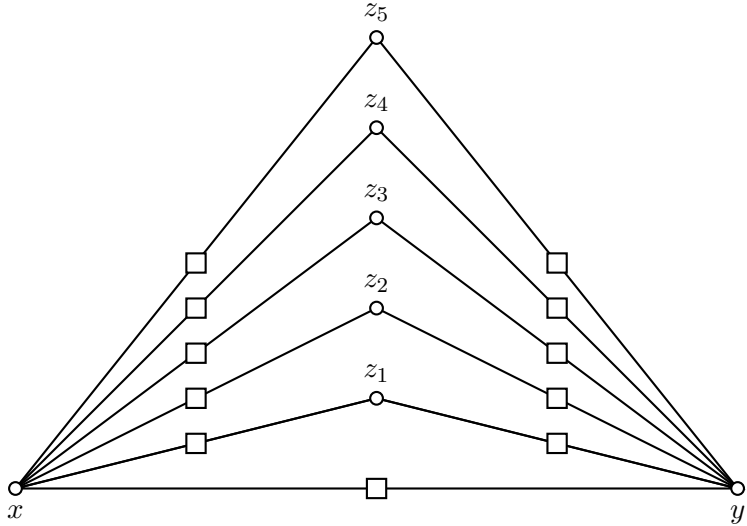
\begin{figure}[ht!]
\begin{center}
\begin{tikzpicture}[scale=0.8,style=thick]
\tikzstyle{every node}=[draw=none,fill=none]
\def\vr{3pt} 

\begin{scope}[yshift = 0cm, xshift = 0cm]
\path (0,0) coordinate (x);
\path (12,0) coordinate (y);
\path (6,1.5) coordinate (z1);
\path (6,3) coordinate (z2);
\path (6,4.5) coordinate (z3);
\path (6,6) coordinate (z4);
\path (6,7.5) coordinate (z5);
\path (6,0) coordinate (xy);
\path (3,0.75) coordinate (xz1);
\path (9,0.75) coordinate (yz1);
\path (3,1.5) coordinate (xz2);
\path (9,1.5) coordinate (yz2);
\path (3,2.25) coordinate (xz3);
\path (9,2.25) coordinate (yz3);
\path (3,3.00) coordinate (xz4);
\path (9,3.00) coordinate (yz4);
\path (3,3.75) coordinate (xz5);
\path (9,3.75) coordinate (yz5);
\draw (x) -- (y) -- (z1) -- (x);
\draw (x) -- (z1) -- (y);
\draw (x) -- (z2) -- (y);
\draw (x) -- (z3) -- (y);
\draw (x) -- (z4) -- (y);
\draw (x) -- (z5) -- (y);
\draw (x)  [fill=white] circle (\vr);
\draw (y)  [fill=white] circle (\vr);
\draw (z1)  [fill=white] circle (\vr);
\draw (z2)  [fill=white] circle (\vr);
\draw (z3)  [fill=white] circle (\vr);
\draw (z4)  [fill=white] circle (\vr);
\draw (z5)  [fill=white] circle (\vr);
\node at (xy) [rectangle,draw,fill=white] () {};
\node at (xz1) [rectangle,draw,fill=white] () {};
\node at (yz1) [rectangle,draw,fill=white] () {};
\node at (xz2) [rectangle,draw,fill=white] () {};
\node at (yz2) [rectangle,draw,fill=white] () {};
\node at (xz3) [rectangle,draw,fill=white] () {};
\node at (yz3) [rectangle,draw,fill=white] () {};
\node at (xz4) [rectangle,draw,fill=white] () {};
\node at (yz4) [rectangle,draw,fill=white] () {};
\node at (xz5) [rectangle,draw,fill=white] () {};
\node at (yz5) [rectangle,draw,fill=white] () {};
\draw[below] (x)++(0.0,-0.1) node {$x$};
\draw[below] (y)++(0.0,-0.1) node {$y$};
\draw[above] (z1)++(0.0,0.1) node {$z_1$};
\draw[above] (z2)++(0.0,0.1) node {$z_2$};
\draw[above] (z3)++(0.0,0.1) node {$z_3$};
\draw[above] (z4)++(0.0,0.1) node {$z_4$};
\draw[above] (z5)++(0.0,0.1) node {$z_5$};
\end{scope}
\end{tikzpicture}
\end{center}
\caption{The graph $S(G_5)$}
\label{fig:difference-2}
\end{figure}

Let $G_n$, $n\ge 2$, be graphs defined as follows. The vertex set of $G_n$ is $\{x,y, z_1, z_2, \ldots, z_n\}$ and the edge set is $\{xy, xz_1, yz_1, xz_2, yz_2,\ldots, xz_n, yz_n\}$. In Fig.~\ref{fig:difference-2} the subdivision $S(G_5)$ of $G_5$ is drawn, where the subdivision vertices are drawn by squares.

One can check that $\mdim(G_2)=4$ and $\mdim(S(G_2))=3$ hold. For instance, $\{v_{xz_1}, v_{xz_2}, v_{yz_1}\}$ is a mixed metric basis of $S(G_2)$. This particular example shows that the right-hand side inequality of Corollary~\ref{cor:subdivided} can be strict. More generally, we have the following result.

\begin{prop}
\label{prop:G_n}
If $n\ge 5$, then $\mdim(G_n) - \mdim(S(G_n)) \ge 2$.
\end{prop}

\proof
Note first that each vertex of $G_n$ has a maximal neighbor, hence by Lemma~\ref{lem:max-neigbor}, $\mdim(G_n) = n+2$. To prove the result it thus suffices to show that $\mdim(S(G_n)) \le n$. For this sake let $S_n = \{v_{xz_1}, v_{xz_2}, v_{yz_3}, v_{yz_4}, z_5, \ldots, z_n\}$. We claim that $S_n$ is a mixed resolving set of $G_n$ and proceed by induction.

For $n=5$ it was checked by hand that the set $S_5$ is indeed a mixed resolving set of $S(G_5)$. Suppose now that the assertion holds for $n-1\ge 5$ and consider $G_n$ and the set $S_n$. Let $x$ and $y$ be two  arbitrary elements from $V(G_n)\cup E(G_n)$. If $x,y$ both lie in the subgraph $G_{n-1}$ of $G_n$, then because $G_{n-1}$ is an isometric subgraph of $G_n$, the elements $x$ and $y$ are by induction distinguished   by some element from $S_{n-1} \subset S_n$. Suppose next that exactly one of $x$ and $y$, say $x$, lies in $G_{n-1}$. Then note that $d_{G_n}(y,z_n) \le 1$ while $d_{G_n}(x,z_n) \ge 2$. Hence $x$ and $y$ are identified by $z_n\in S_n$. The last case to consider is when none of $x$ and $y$ lies in $G_{n-1}$. In that case, if $d_{G_n}(x,z_n) = d_{G_n}(y,z_n)$, then $d_{G_n}(x,v_{xz_1}) \ne d_{G_n}(y,v_{xz_1})$ and we are done.
\qed

We would point out an open problem concerning whether there are graphs $G$ such that $\mdim(G) - \mdim(S(G)) > 2$. On the other hand, equality cases are also possible, for instance we have the following:

\begin{cor}
If $G$ is a cactus graph, then $\mdim(G) = \mdim(S(G))$.
\end{cor}

\proof
It is straightforward to verify that $S(G)$ is a cactus graph such that (i) $C$ is a cycle in $G$ if and only if $S(C)$ is a cycle in $S(G)$; (ii)  $\rt(C)\geq3$ and there is not a geodesic triple of root vertices on the cycle $C$ if and only if $\rt(S(C))\geq 3$ and there is not a geodesic triple of root vertices on the cycle $S(C)$; (iii) and that $\rt(C)=\rt(S(C))$. Therefore,  Theorem~\ref{th0} implies the result.
\qed

Another interesting open question is to characterize graphs $G$ such that $\mdim(G) = \mdim(S(G))$.

\section{Middle and total graphs}

In the first part of this section we deal with middle graphs. We first demonstrate that $\dim(M(G))$ can be bounded from the above by $\mdim(G)$ and then prove that the bound is tight on trees.

\begin{thm}\label{nco1}
If $G$ is a graph, then $\mdim(G)\geq\dim(M(G))$.
\end{thm}

\proof
We first consider the distances between elements of $S(G)$. If $x,y\in V(G)$, then clearly
\begin{equation}\label{eq:Mid-vertex-vertex}
d_{M(G)}(x,y)=d_G(x,y)+1.
\end{equation}
On the other hand, if $x\in V(G)$ and $e\in E(G)$, then
\begin{equation}\label{eq:Mid-vertex-edge}
d_{M(G)}(x,v_e)=d_G(x,e)+1.
\end{equation}

Next, let $W$ be a mixed resolving set for $G$. We claim that $M$ forms a resolving set for $M(G)$. To this end, let $x,y\in V(M(G))$. The cases $x\in W$ or $y\in W$ are straightforward. Hence, we assume $x,y\notin W$. Since each of $x,y$ corresponds either to a vertex or to an edge of $G$, and $W$ is a mixed resolving set for $G$, there exists a vertex $w\in W$ such that $d_G(x,w)\ne d_G(y,w)$.

If $x,y\in V(G)$, then by \eqref{eq:Mid-vertex-vertex} we have
\[d_{M(G)}(x,w)=d_G(x,w)+1\ne d_G(y,w)+1=d_{M(G)}(y,w)\]
 If $x=v_{e}$ and $y=v_f$ for some $e,f\in E(G)$, then by \eqref{eq:Mid-vertex-edge},
\[d_{M(G)}(v_e,w)=d_G(e,w)+1\ne d_G(f,w)+1=d_{M(G)}(v_f,w).\]
 Finally, if $x\in V(G)$ and $y=v_f$ for some $f\in E(G)$, then by \eqref{eq:Mid-vertex-vertex} and \eqref{eq:Mid-vertex-edge}, it follows
\[d_{M(G)}(x,w)=d_G(x,w)+1\ne d_G(f,w)+1=d_{M(G)}(v_f,w).\]
 Therefore, $W$ is a resolving set for $M(G)$ as claimed.
\qed

\begin{thm}\label{nth1}
If $G$ is a tree, then $\mdim(G)=\dim(M(G))=n_1(G)$.
\end{thm}

\proof
From \cite{KKTY-2017} we have $\mdim(G)=n_1(G)$. Thus, by Theorem~\ref{nco1}, $\dim(M(G))\leq n_1(G)$. To complete our proof it is enough to show that $\dim(M(G))\geq n_1(G)$. To do this, let $W$ be a metric basis for $M(G)$.

Consider a pendant edge $xy$ of $G$, where $d_{G}(x) = 1$. Then $d_{M(G)}(v_{xy}, z) = d_{M(G)}(y, z)$ holds for every vertex $z\in V(M(G))\setminus \{x,v_{xy},y\}$. It follows that $W\cap \{x,v_{xy},y\} \ne \emptyset$. Suppose now that $y$ is a support vertex of $G$ with $k\ge 2$ leaves attached to it, say $x_1, \ldots, x_k$. By the above, $W\cap \{x_i,v_{x_iy},y\} \ne \emptyset$ and $W\cap \{x_j,v_{x_jy},y\} \ne \emptyset$ hold for each $i,j\in [k]$, $i\ne j$. If $W\cap \{x_i, v_{x_iy}, y\} = \{y\} = W\cap \{x_j, v_{x_jy}, y\}$, then $x_i$ and $x_j$ are not distinguished by $W$. Therefore we deduce that
$$|W \cap \bigcup_{i=1}^k \{x_i, v_{x_iy}, y\}| \ge k\,.$$
Hence $W$ contains at least $n_1$ vertices and we are done.
\qed

We now turn our attention to the mixed metric dimension of total graphs of trees and prove the following somewhat surprising result.

\begin{thm}\label{thm:total}
If $G$ is a tree, then
$$\mdim(T(G)) = 2 n_1(G)\,.$$
\end{thm}

\proof
Let $x$ be a leaf of $G$ adjacent to $y$. Then we claim that both $x$ and $v_{xy}$ have a maximal neighbor. Indeed, $v_{xy}$ is a maximal neighbor of $x$, and $y$ is a maximal neighbor of $v_{xy}$. By Lemma~\ref{lem:max-neigbor}, both $x$ and $v_{xy}$ belong to each mixed resolving set. Applying this argument to each pendant edge of $G$ we obtain that $\mdim(T(G)) \ge 2n_1(G)$.

Let now $W$ be the set of vertices of $T(G)$ containing all the vertices corresponding to the leaves in $G$ together with the subdivision vertices adjacent to them. Then $|W| = 2n_1(G)$. We claim that $W$ is a mixed metric basis of $T(G)$. Let $x$ and $y$ be arbitrary elements from $V(T(G))\cup E(T(G))$.

Consider first all possible situations when $x,y\in V(T(G))$. If $x,y\in V(G)$, then there exists a leaf $w$ in $G$ such that $d_G(x,w) \ne d_G(y,w)$. Since $d_{T(G)}(x,w) = d_G(x,w)$ and $d_{T(G)}(y,w) = d_G(y,w)$, $x$ and $y$ are distinguished by $w\in W$ in $T(G)$ also. Suppose next that $x\in V(G)$ and $y = v_{zz'}$. Let $t$ be a leaf of $G$ such that $d_G(x,t)\le \min\{d_G(z,t), d_G(z',t)\}$. Then $d_{T(G)}(y,t) \ge d_{T(G)}(x,t) + 1$. The next case is when $x=v_{zz'}$ and $y=v_{ww'}$, which can be dealt with using the same argument as in the previous case.

Consider next all possible situations when $x,y\in E(T(G))$. If $x$ and $y$ are original edges, then their distances to leaves remain the same and there is nothing to prove, based on the facts that the sets of leaves in $G$ and $T(G)$ are the same, that the set of leaves of $G$ is included in $W$, and that the set of leaves is a mixed metric basis of $G$.

We next consider all the remaining cases when $x$ and $y$ are incident edges.

\medskip\noindent
{\bf Case 1}: $x$ and $y$ are $S(G)$-edges: $x=zv_{zz'}$ and $y=z'v_{zz'}$.\\
Then let $t$ be a leaf of $G$ such that $d_G(t,z') = d_G(t,z) + 1$. Then
\begin{equation*}
\label{eq:to-leaf}
d_{T(G)}(y,t) = d_{T(G)}(x,t) + 1
\end{equation*}
which means that $x$ and $y$ are distinguished by $t$ which is a leaf of $G$.

\medskip\noindent
{\bf Case 2}: $x$ and $y$ are $S(G)$-edges: $x = z'v_{zz'}$ and $y = z'v_{z'z''}$. \\
Then select a leaf $t$ of $G$ such that the shortest path from $t$ to $z''$ contains $z$ and $z'$. Consider now the vertex $v_{tt'}\in W$, where $tt'\in E(G)$. Then we infer that
\begin{equation*}
\label{eq:to-neighbor-of-leaf}
d_{T(G)}(y,v_{tt'}) = d_{T(G)}(x,v_{tt'}) + 1
\end{equation*}
which means that $x$ and $y$ are distinguished by the subdivision vertex $v_{tt'}$ which is adjacent to the leaf $t$ of $G$.

\medskip
For all the other situations in which $x$ and $y$ are incident, we consider Cases 1 and 2 as models. We have collected the remaining situations in Table~\ref{tab:all-cases}, where the column corresponding to the identifying vertex contains a vertex from $W$ selected in a way parallel as the vertices $t$ and $v_{tt'}$ were selected in Cases 1 and 2, respectively.

\begin{table}[ht!]
\centering
\begin{tabular}{|c| c | c |}
 \hline
edge $x$ of type & edge $y$ of type & identifying vertex \\
[0.5ex]
\hline\hline
$L(G): v_{zz'}v_{z'z''}$ & $L(G): v_{z'z''}v_{z''z'''}$ & \phantom{$X^{X^X}$} $v_{tt'}$ \phantom{$X^{X^X}$} \\
[0.5ex]
\hline
$L(G): v_{zz'}v_{z'z''}$ & $S(G): v_{z'z''}z''$ & \phantom{$X^{X^X}$} $v_{tt'}$ \phantom{$X^{X^X}$} \\
[0.5ex]
\hline
$L(G): v_{zz'}v_{z'z''}$ & $S(G): v_{z'z''}z'$ & \phantom{$X^{X^X}$} $v_{tt'}$ \phantom{$X^{X^X}$} \\
[0.5ex]
\hline
original: $zz'$ & $S(G): v_{zz'}z'$ & \phantom{$X^{X^X}$} $t$ \phantom{$X^{X^X}$} \\
[0.5ex]
\hline
original: $zz'$ & $S(G): z'v_{z'z''}$ & \phantom{$X^{X^X}$} $t$ \phantom{$X^{X^X}$} \\
[0.5ex]
\hline
\end{tabular}
\caption{Some remaining pairs of edges and one identifying vertex of each pair.}
\label{tab:all-cases}
\end{table}

To cover all the cases when $x$ and $y$ are edges, we still need to deal with the situation when $x$ and $y$ are not incident. For instance, this is the case when $x$ is an original edge and $y$ is an $L(G)$-edge. However, for non-incident edges the existence of a vertex of type $t$ or  $v_{tt'}$ can be readily deduced.

To complete the proof, we need to see that also an edge $x$ and a vertex $y$ are distinguished by $W$. Suppose first that $x$ is a vertex of $G$ and $y = xv_{zx}$. Then we consider a leaf $t$ of $G$ and the subdivision vertex $v_{tt'}\in W$, where $tt'\in E(G)$, such that $d_{T(G)}(x,v_{tt'}) = d_{T(G)}(y,v_{tt'}) + 1$ and we are done. Suppose second that $x = v_{zz'}$ and $y=zv_{zz'}$. Then we detect a leaf $t$ such that $d_{T(G)}(x,t) = d_{T(G)}(y,t) + 1$. In any other case (including those in which $x$ and $y$ are not incident) we can similarly find a vertex $t$ (leaf of $G$) or $v_{tt'}$ (subdivision vertex of $T(G)$ adjacent to the leaf $t$ of $G$) which distinguish $x$ and $y$.
\qed

To prove the left-hand side inequality of our next result, we use the standard terminology concerning the metric dimension of trees that can be seen in, say~\cite{HM-1976, dorota-2022+, S-1975, till-2022+}. That is, a vertex of degree at least three of a tree $G$ is called a \emph{major vertex}. A leaf $u$ of $G$ is called a \emph{terminal vertex} of a major vertex $v$ of $G$ if $d_G(u,v) < d_G(u,w)$ for every other major vertex $w$ of $G$. The \emph{terminal degree} of a major vertex $v$ is the number of terminal vertices of $v$. A major vertex $v$ of $G$ is an \emph{exterior major vertex} of $G$ if it has positive terminal degree.

\begin{prop}
\label{prop:dim-T(G)}
If $G$ is a tree, then $\dim(G) \le \dim(T(G)) \le n_1(G)$.
\end{prop}

\proof
In the third paragraph of the proof of Theorem~\ref{thm:total} when considering arbitrary two vertices of $T(G)$, each pair of vertices was distinguished by some leaf of $G$. This implies the right-hand side inequality.

Now, for any exterior major vertex $u$ of $G$, the subtree $G_u$ containing all the vertices between $u$ and its terminal leaves, must contain exactly the number of terminal leaves of $u$ minus $1$. Denote this number by $g_u$.  To prove our assertion, we claim that a basis $W$ of $T(G)$ must contain at least $g_u$ vertices from $T(G_u)$. Suppose this is not the case. Then there exist two leaves $t_1$ and $t_2$ in $G$ such that no vertex in $T(G)$ on any shortest $t_1,u$-path and no vertex on any shortest $t_2,u$-path  lies in $W$. But then the vertices $v_{uu'}$ and $v_{uu''}$ are not identified by $W$, where $u'$ and $u''$ are adjacent to $u$ in $G$ on the $u,t_1$- and $u,t_2$-paths, respectively.
\qed

The upper bound of Proposition~\ref{prop:dim-T(G)} is tight as demonstrated by paths and stars $K_{1,k}$, $k\in \{2,3,4\}$ for which $\dim(T(P_n)) = 2$ and $\dim(T(K_{1,k})) = k$ for $k\in \{2,3,4\}$, see~\cite{soo-2016}. On the other hand, $\dim(T(K_{1,k})) = k-1$ for $k\ge 5$, again see~\cite{soo-2016}, which shows that the lower bound is also tight.  Hence the problem to determine $\dim(T(G))$ for an an arbitrary tree seems a challenging problem, an interesting fact because by Theorem~\ref{thm:total} we know $\mdim(T(G))$.


\section*{Acknowledgments}

Sandi Klav\v{z}ar acknowledges the financial support from the Slovenian Research Agency through research core funding No.\ P1-0297 and projects J1-2452 and N1-0285. Ismael G. Yero has been partially supported by the Spanish Ministry of Science and Innovation through the grant PID2019-105824GB-I00. Moreover, this investigation was developed while this last author was visiting the University of Ljubljana, Slovenia, supported by ``Ministerio de Educaci\'on, Cultura y Deporte'', Spain, under the ``Jos\'e Castillejo'' program for young researchers (reference number: CAS21/00100).



\begin{thebibliography}{99}

\bibitem{behzad-1966}
M.~Behzad, G.~Chartrand,
Total graphs and traversability,
Proc.\ Edinburgh Math.\ Soc.\ 15 (1966/67) 117--120.

\bibitem{FKK-2019}
V.~Filipovi\'c, A.~Kartelj, J.~Kratica,
Edge metric dimension of some generalized Petersen graphs,
Results Math.\ 74 (2019) Paper No.\ 182.

\bibitem{hamada-1976}
T.~Hamada, I.~Yoshimura,
Traversability and connectivity of the middle graph of a graph,
Discrete Math.\ 14 (1976) 247--255.

\bibitem{HM-1976}
F. Harary, R.A.~Melter,
On the metric dimension of a graph,
Ars Combin.\ 2 (1976) 191--195.

\bibitem{Johnson-1993} M.~Johnson,
Structure-activity maps for visualizing the graph variables arising in drug design,
J.\ Biopharm.\ Statist.\ 3(2) (1993) 203--236.

\bibitem{KTY-2018}
A.~Kelenc, N.~Tratnik, I.G.~Yero,
Uniquely identifying the edges of a graph: the edge metric dimension,
Discrete Appl.\ Math.\ 251 (2018) 204--220.

\bibitem{KKTY-2017}
A.~Kelenc, D.~Kuziak, A.~Taranenko, I.G.~Yero,
Mixed metric dimension of graphs,
Appl.\ Math.\ Comput.\ 314 (2017) 429--438.

\bibitem{dorota-2022+}
D.~Kuziak, I.G.~Yero,
Metric dimension related parameters in graphs: A survey on combinatorial, computational and applied results,
arXiv:2107.04877 [math.CO] (10 Jul 2021).

\bibitem{kwak-2006}
J.H.~Kwak, I.~Sato,
Zeta functions of line, middle, total graphs of a graph and their coverings,
Linear Algebra Appl.\ 418 (2006) 234--256.

\bibitem{MM2}
J.B.~Liu, S.~Kumar Sharma, V.~Kumar Bhat, H.~Raza,
Multiset and mixed metric dimension for starphene and zigzag-edge coronoid,
Polycyclic Arom.\ Comp.\ (2021) doi.org/10.1080/10406638.2021.2019066.

\bibitem{Melter-1984}
R.A.~Melter, I.~Tomescu,
Metric bases in digital geometry,
Comput.\ Vision Graphics Image Process.\ 25(1) (1984) 113--121.

\bibitem{MKSM-2021}
M.~Milivojevi\'c Danas, J.~Kratica, A.~Savi\'c, Z.Lj.~Maksimovi\'c,
Some new general lower bounds for mixed metric dimension of graphs,
Filomat 35 (2021) 4275--4285.

\bibitem{SS-2021}
J.~Sedlar, R.~\v{S}krekovski,
Mixed metric dimension of graphs with edge disjoint cycles,
Discrete Appl.\ Math.\ 300 (2021) 1--8.

\bibitem{sedlar-2022}
J.~Sedlar, R.~\v{S}krekovski,
Vertex and edge metric dimensions of unicyclic graphs,
Discrete Appl.\ Math.\ 314 (2022) 81--92.

\bibitem{S-1975}
P.J.~Slater,
Leaves of trees,
Cong.\ Numer.\ 14 (1975) 549--559.

\bibitem{somodi-2017}
M.~Somodi, K.~Burke, J.~Todd,
On a construction using commuting regular graphs,
Discrete Math.\ 340 (2017) 532--540.

\bibitem{soo-2016}
B.~Sooryanarayana, Shreedhar~K., Narahari~N.,
On the metric dimension of the total graph of a graph,
Notes Number Theory Discrete Math.\ 22(4) (2016) 82--95.

\bibitem{till-2022+}
R.C.~Tillquist, R.M.~Frongillo, M.E.~Lladser,
Getting the lay of the land in discrete space: A survey of metric dimension and its applications,
arXiv:2104.07201 [math.CO] (15 Apr 2021).

\bibitem{till-2019}
R.C.~Tillquist, M.E.~Lladser,
Low-dimensional representation of genomic sequences,
J.\ Math.\ Biol.\ 79(1) (2019) 1--29.

\bibitem{Trujillo-Rasua-2016}
R.~Trujillo-Ras\'ua, I.G.~Yero,
$k$-metric antidimension: A privacy measure for social graphs,
Inform.\ Sci.\ 328 (2016) 403--417.


\end{thebibliography}
\end{document}